\documentclass[12pt]{article}

\usepackage{amsmath,amsfonts,amssymb,latexsym,amscd}

\begin{document}

\newtheorem{theorem}{Theorem}[section]
\newtheorem{corollary}[theorem]{Corollary}
\newtheorem{lemma}[theorem]{Lemma}
\newtheorem{proposition}[theorem]{Proposition}
\newtheorem{conjecture}[theorem]{Conjecture}
\newtheorem{commento}[theorem]{Comment}
\newtheorem{definition}[theorem]{Definition}
\newtheorem{problem}[theorem]{Problem}
\newtheorem{remark}[theorem]{Remark}
\newtheorem{remarks}[theorem]{Remarks}
\newtheorem{example}[theorem]{Example}

\newcommand{\Nb}{{\mathbb{N}}}
\newcommand{\Rb}{{\mathbb{R}}}
\newcommand{\Tb}{{\mathbb{T}}}
\newcommand{\Zb}{{\mathbb{Z}}}
\newcommand{\Cb}{{\mathbb{C}}}

\newcommand{\Ef}{\mathfrak E}
\newcommand{\Gf}{\mathfrak G}
\newcommand{\iGf}{\mathfrak I\mathfrak G}
\newcommand{\Hf}{\mathfrak H}
\newcommand{\Kf}{\mathfrak K}
\newcommand{\Lf}{\mathfrak L}
\newcommand{\Af}{\mathfrak A}
\newcommand{\Bf}{\mathfrak B}

\def\A{{\mathcal A}}
\def\B{{\mathcal B}}
\def\C{{\mathcal C}}
\def\D{{\mathcal D}}
\def\F{{\mathcal F}}
\def\G{{\mathcal G}}
\def\H{{\mathcal H}}
\def\J{{\mathcal J}}
\def\K{{\mathcal K}}
\def\LL{{\mathcal L}}
\def\N{{\mathcal N}}
\def\M{{\mathcal M}}
\def\N{{\mathcal N}}
\def\O{{\mathcal O}}
\def\P{{\mathcal P}}
\def\SS{{\mathcal S}}
\def\T{{\mathcal T}}
\def\U{{\mathcal U}}
\def\W{{\mathcal W}}

\def\span{\operatorname{span}}
\def\Ad{\operatorname{Ad}}
\def\ad{\operatorname{Ad}}
\def\tr{\operatorname{tr}}
\def\id{\operatorname{id}}
\def\en{\operatorname{End}}
\def\aut{\operatorname{Aut}}
\def\out{\operatorname{Out}}
\def\per{\operatorname{Per}(X_n)}
\def\la{\langle}
\def\ra{\rangle}

\def\j{j_\infty}
\def\f{f_\infty}
\def\g{g_\infty}
\def\a{a_\infty}

\title{The Restricted Weyl Group of the Cuntz Algebra and Shift Endomorphisms}

\author{Roberto Conti, Jeong Hee Hong\footnote{This work was supported by 
National Research Foundation of Korea
Grant funded by the Korean Government (KRF--2008--313-C00039).} and
Wojciech Szyma{\'n}ski\footnote{This work was partially supported by the FNU Rammebevilling
grant `Operator algebras and applications' (2009--2011),  the NordForsk Research 
Network `Operator algebra and dynamics' (2009--2011), and the Marie Curie Research 
Training Network MRTN-CT-2006-031962 EU-NCG.}}

\date{{\small 22 June 2010}}
\maketitle

\renewcommand{\sectionmark}[1]{}

\vspace{7mm}
\begin{abstract}
It is shown that, modulo the automorphisms which fix the canonical diagonal MASA point-wise, 
the group of those automorphisms of $\O_n$ which globally preserve both the diagonal and 
the core UHF-subalgebra is isomorphic, via restriction, with the group of those homeomorphisms 
of the full one-sided $n$-shift space which eventually commute along with their inverses with 
the shift transformation. The image of this group in the outer automorphism group of $\O_n$ 
can be embedded into the quotient of the automorphism group of the full two-sided 
$n$-shift by its center, generated by the shift. If $n$ is prime then this embedding is an 
isomorphism. 
\end{abstract}

\vfill\noindent {\bf MSC 2010}: 46L40, 37B10

\vspace{3mm}

\newpage

\section{Introduction}

Investigations of endomorphisms of the simple $C^*$-algebras $\O_n$, \cite{Cun1}, were initiated 
by Cuntz in his seminal paper \cite{Cun2}. A key tool for those investigations was provided by a fairly simple 
observation that unital endomorphisms of $\O_n$ are in an explicit bijective correspondence with 
unitary elements of this $C^*$-algebra. In particular, Cuntz studied the subgroup of automorphisms 
of $\O_n$ which globally preserve the canonical diagonal MASA $\D_n$. He showed that the quotient 
of this group by its normal subgroup consisting of those automorphisms which fix $\D_n$ point-wise (and 
this is a maximal abelian subgroup of $\aut(\O_n)$) is discrete. It is natural to think of this quotient as an 
analogue of the Weyl group. The Weyl group contains a natural interesting subgroup corresponding to 
those automorphisms which also globally preserve the core UHF-subalgebra $\F_n$ of $\O_n$. Cuntz 
proposed in \cite{Cun2} a problem of determining the structure of this restricted Weyl group. 

Ever since, endomorphisms of $\O_n$ have been an active area of investigations both intrinsically and 
in connection to many other areas. Not even attempting any exhausting overview of the relevant 
literature, let us only mention two pieces of research which influenced our current study the most: 
Conti-Pinzari work on Jones index for endomorphisms of $\O_n$, \cite{CP}, and Bratteli-J{\o}rgensen 
work on iterated function systems and representations of Cuntz algebras, initiated in \cite{BJ1} and 
\cite{BJ2}. More recently, quite significant progress has been achieved in the study of those endomorphisms 
which preserve either the core UHF-subalgebra, \cite{CRS}, or the diagonal MASA, \cite{HSS}. A powerful, 
novel combinatorial approach to the study of endomorphisms which globally preserve both $\F_n$ and 
$\D_n$ has been developed in \cite{S}, \cite{CS} and \cite{CKS}. 

The main result of this paper is an explicit and intrinsic description 
of the restricted Weyl group of the Cuntz algebra $\O_n$ and its image in 
the outer automorphism group $\out(\O_n)$ (which we call the restricted outer Weyl group of $\O_n$). 
The way we achieve this is by 
analyzing the action of the restricted Weyl group on the diagonal MASA, and by showing 
that certain class of automorphisms of $\D_n$ admits extensions to permutative automorphisms 
of $\O_n$ (and thus also of $\F_n$). In general, existence of such extensions is 
not guaranteed, as demonstrated by \cite{Co}. 
It turns out that the restricted outer Weyl group admits a natural embedding into 
the group of shift automorphisms of the two-sided full shift (with the mod out center), 
and for prime $n$ this embedding is actually an isomorphism. These facts 
have profound implications. For example, they immediately imply that the 
restricted outer Weyl group is residually finite. 
For $n =2$, this result also provides a not unexpected answer to a 
question left open in \cite{CS}. Namely, the restricted outer Weyl 
group is not amenable (for $n \geq 3$ this was already shown earlier in \cite{S} and \cite{CS}).
It also follows from our results that the restricted Weyl group of $\O_n$ is big 
enough to contain a copy of the group of shift automorphisms of the 
one-sided full $n$-shift. The groups of shift automorphisms (both for the one-sided and 
the two-sided shift) have been extensively studied in the literature as they reveal an intriguing and 
highly nontrivial structure, see \cite{K,LM} and the literature cited therein.

One important outcome of our investigations is a powerful and neat link between two areas: 
the study of automorphisms of the Cuntz algebras and symbolic dynamics. This kind of interaction 
was pioneered by Cuntz and Krieger in \cite{CK}, but we believe its limits have not been 
reached yet. On one hand, we expect that the nice results available in the literature 
on symbolic dynamics could shed new light on some aspects of automorphisms of Cuntz 
algebras, and possibly even more general classes of $C^*$-algebras.
But perhaps it is even more intriguing to speculate the other way round, as one could hope that 
the algebraic environment that we unveil could provide new tools 
and insight for attacking some of the problems in the dynamical systems setting.

\vspace{3mm}\noindent
{\bf Acknowledgements.} 
We are utmost grateful to Mike Boyle for a very helpful and illuminating discussion of shift endomorphisms.

\section{Notation and preliminaries}

If $n$ is an integer greater than 1, then the Cuntz algebra $\O_n$ is a unital, simple,
purely infinite $C^*$-algebra generated by $n$ isometries $S_1, \ldots, S_n$, satisfying
$\sum_{i=1}^n S_i S_i^* = I$, \cite{Cun1}.
We denote by $W_n^k$ the set of $k$-tuples $\mu = (\mu_1,\ldots,\mu_k)$
with $\mu_m \in \{1,\ldots,n\}$, and by $W_n$ the union $\cup_{k=0}^\infty W_n^k$,
where $W_n^0 = \{0\}$. We call elements of $W_n$ multi-indices.
If $\mu \in W_n^k$ then $|\mu| = k$ is the length of $\mu$.
If $\mu = (\mu_1,\ldots,\mu_k) \in W_n$, then $S_\mu = S_{\mu_1} \ldots S_{\mu_k}$
($S_0 = 1$ by convention) is an isometry with range projection $P_\mu=S_\mu S_\mu^*$.
Every word in $\{S_i, S_i^* \ | \ i = 1,\ldots,n\}$ can be uniquely expressed as
$S_\mu S_\nu^*$, for $\mu, \nu \in W_n$ \cite[Lemma 1.3]{Cun1}.

We denote by $\F_n^k$ the $C^*$-subalgebra of $\O_n$ spanned by all words of the form
$S_\mu S_\nu^*$, $\mu, \nu \in W_n^k$, which is isomorphic to the
matrix algebra $M_{n^k}({\mathbb C})$. The norm closure $\F_n$ of
$\cup_{k=0}^\infty \F_n^k$, is the UHF-algebra of type $n^\infty$,
called the core UHF-subalgebra of $\O_n$, \cite{Cun1}. It is the fixed point algebra
for the gauge action of the circle group $\gamma:U(1)\rightarrow{\aut}(\O_n)$ defined
on generators as $\gamma_t(S_i)=tS_i$. For $k\in\Zb$, we denote by $\O_n^{(k)}
:=\{x\in\O_n:\gamma_t(x)=t^k x\}$,
the spectral subspace for this action. In particular, $\F_n=\O_n^{(0)}$.
The $C^*$-subalgebra of $\F_n$ generated by projections $P_\mu$, $\mu\in W_n$, is a 
MASA (maximal abelian subalgebra) both in $\F_n$ and in $\O_n$. We call it
the diagonal and denote $\D_n$. The spectrum of $\D_n$ is naturally identified with 
$X_n$ --- the full one-sided $n$-shift space. We also set $\D_n^k:=\D_n\cap\F_n^k$. 
Throughout this paper we are interested in the inclusions
$$ \D_n \subseteq \F_n\subseteq \O_n.$$
The UHF-subalgebra $\F_n$ posseses a unique normalized trace, denoted $\tau$. 
We will refer to the restriction of $\tau$ to $\D_n$ as to the canonical trace on $\D_n$. 

We denote by $\SS_n$ the group of those unitaries in $\O_n$ which can be written
as finite sums of words, i.e., in the form $u = \sum_{j=1}^m S_{\mu_j}S_{\nu_j}^*$
for some $\mu_j, \nu_j \in W_n$. We also denote $\P_n=\SS_n\cap\U(\F_n)$. Then
$\P_n=\cup_k\P_n^k$, where $\P_n^k$ are permutation unitaries in $\U(\F_n^k)$.
That is, for each $u\in\P_n^k$ there is a unique permutation $\sigma$ of multi-indices
$W_n^k$ such that
\begin{equation}\label{permutunitary}
u = \sum_{\mu \in W_n^k} S_{\sigma(\mu)} S_\mu^*.
\end{equation}

As shown by Cuntz in \cite{Cun2}, there exists the following bijective correspondence
between unitaries in $\O_n$ and unital $*$-endomorphisms of $\O_n$ (whose collection we denote
by $\en(\O_n)$). A unitary $u$  in $\O_n$ determines an endomorphism $\lambda_u$ by
$$ \lambda_u(S_i) = u S_i, \;\;\; i=1,\ldots, n. $$
Conversely, if $\rho :\O_n\rightarrow \O_n$ is an endomorphism, then
$\sum_{i=1}^n\rho(S_i)S_i^*=u$ gives a unitary $u\in\O_n$
such that $\rho=\lambda_u$. If the unitary $u$ arises from a permutation $\sigma$ via the formula
\eqref{permutunitary}, the corresponding endomorphism will be 
sometimes denoted by $\lambda_{\sigma}$.
Composition of endomorphisms corresponds to a `convolution'
multiplication of unitaries as follows:
\begin{equation}\label{convolution}
\lambda_u \circ \lambda_w = \lambda_{\lambda_u(w)u}
\end{equation}
We denote by $\varphi$ the canonical shift:
$$ \varphi(x)=\sum_iS_ixS_i^*, \;\;\; x\in\O_n. $$
If we take $u=\sum_{i, j}S_iS_jS_i^*S_j^*$ then $\varphi=\lambda_u$.
It is well-known that $\varphi$ leaves invariant both $\F_n$ and $\D_n$, and that $\varphi$
commutes with the gauge action $\gamma$. We denote by $\phi$ the standard left inverse 
of $\varphi$, defined as $\phi(a)=\frac{1}{n}\sum_{i=1}^n S_i^* a S_i$. 

If $u\in\U(\O_n)$ then for each positive integer $k$ we denote
\begin{equation}\label{uk}
u_k = u \varphi(u) \cdots \varphi^{k-1}(u).
\end{equation}
We agree that $u_k^*$ stands for $(u_k)^*$. If
$\alpha$ and $\beta$ are multi-indices of length $k$ and $m$, respectively, then
$\lambda_u(S_\alpha S_\beta^*)=u_kS_\alpha S_\beta^*u_m^*$. This is established through
a repeated application of the identity $S_i a = \varphi(a)S_i$, valid for all
$i=1,\ldots,n$ and $a \in \O_n$. If $u\in\F_n^k$ for some $k$ then, 
following \cite{CP},  we call endomorphism $\lambda_u$ {\em localized}. 

For algebras $A\subseteq B$ we denote by $\N_B(A)=\{u\in\U(B):uAu^*=A\}$ the normalizer
of $A$ in $B$ and by $A' \cap B=\{a \in A: (\forall b \in B) \; ab=ba\}$ the 
relative commutant of $A$ in $B$. We also denote by $\aut(B,A)$ the collection of all those 
automorphisms $\alpha$ of $B$ that $\alpha(A)=A$, and by $\aut_A(B)$ those 
automorphisms of $B$ which fix $A$ point-wise. 

%%%%%%%%%%%%%%%%%%%%%%%%%%%%%%%%%%%%%%%%%%%%%
%%%%%%%%%%%%%%%%%%%%%%%%%%%%%%%%%%%%%%%%%%%%%

\section{The restricted Weyl group of the Cuntz algebra}

Let $\alpha$ be an automorphism of $\D_n$. We denote by $\alpha_*$ the corresponding 
homeomorphism of $X_n$. We say that $\alpha$ has {\em property (P)} if 
there exists $m$ such that for all $k\geq m$ we have 
\begin{equation}\label{propertyP}
\alpha\varphi^k(x)=\varphi^{k-m}\alpha\varphi^m(x)
\end{equation}
for all $x\in\D_n^1$. That is, $\alpha$ eventually commutes with the shift. Equivalently, 
$\alpha\in\aut(\D_n)$ satisfies (P) with $m$ if the endomorphism $\alpha\varphi^m$ commutes 
with the shift $\varphi$. We define 
\begin{equation}\label{gn}
\Gf_n:=\{\alpha\in\aut(\D_n):\text{ both $\alpha$ and $\alpha^{-1}$ have property (P)}\}. 
\end{equation}

\begin{lemma}\label{groupgn}
$\Gf_n$ is a subgroup of $\aut(\D_n)$. 
\end{lemma}
{\em Proof.} 
Let $\alpha$ and $\beta$ belong to $\Gf_n$. Take $m$ so large that both $\alpha$ and $\beta$ 
satisfy (P) with $m$. Also, let $r$ be so large that $\beta(\varphi^m(\D_n^1))$ is contained in 
$\D_n^r$. Since $\D_n^r=\D_n^1\varphi(\D_n^1)\cdots\varphi^{r-1}(\D_n^1)$, there exist 
linear functionals $f_\mu:\D_n^1\to\Cb$, $\mu\in W_n^r$, such that 
\begin{equation}\label{betam}
\beta(\varphi^m(x))=\sum_{\mu\in W_n^r}f_\mu(x)P_{\mu_0}\varphi(P_{\mu_1}) 
\cdots \varphi^{r-1}(P_{\mu_{r-1}}). 
\end{equation}
Thus for $k\geq 2m$ and $x\in\D_n^1$ we have 
\begin{eqnarray*}
\alpha\beta\varphi^k(x) & = & \alpha\varphi^{k-m}\beta\varphi^m(x) \\
 & = & \alpha\varphi^{k-m}\sum_{\mu\in W_n^r} f_\mu(x)P_{\mu_0}\varphi(P_{\mu_1}) 
\cdots \varphi^{r-1}(P_{\mu_{r-1}}) \\
 & = & \sum_{\mu\in W_n^r} f_\mu(x)(\alpha\varphi^{k-m}(P_{\mu_0}))
 \cdots (\alpha\varphi^{k-m+r-1}(P_{\mu_{r-1}})) \\
 & = & \sum_{\mu\in W_n^r} f_\mu(x)(\varphi^{k-2m}\alpha\varphi^{2m}(P_{\mu_0}))
 \cdots (\varphi^{k-2m}\alpha\varphi^{2m+r-1}(P_{\mu_{r-1}})) \\
 & = & \varphi^{k-2m}\alpha\sum_{\mu\in W_n^r} f_\mu(x)\varphi^{2m}(P_{\mu_0})
 \cdots \varphi^{2m+r-1}(P_{\mu_{r-1}}) \\
 & = & \varphi^{k-2m}\alpha\beta\varphi^{2m}(x),  
\end{eqnarray*}
since 
$$ \sum_{\mu\in W_n^r} f_\mu(x)\varphi^{2m}(P_{\mu_0})
 \cdots \varphi^{2m+r-1}(P_{\mu_{r-1}}) =  \beta\varphi^{2m}(x) $$ 
by virtue of (\ref{betam}) and property (P) for $\beta$. Consequently, the product $\alpha\beta$ 
satisfies (P) with $2m$, and whence $\Gf_n$ is a group. 
\hfill $\Box$ 

\medskip
We denote $\iGf_n=\{\ad(u)|_{\D_n}:u\in\P_n\}$. This is a normal subgroup of $\Gf_n$, 
since for $u\in\P_n^k$ we have $\ad(u)\varphi^k=\varphi^k$. In what follows we agree that $\varphi^0=\id$.

\begin{lemma}\label{amk}
If $\alpha\in\aut(\D_n)$ then $\alpha\in\iGf_n$ if and only if there exist $m,k$ such that $\alpha\varphi^m=\varphi^k$.  
\end{lemma}
{\em Proof.} Let $\alpha\in\aut(\D_n)$ and let $m,k$ be such that $\alpha\varphi^m=\varphi^k$. 
Suppose first that $m\geq k$. Then $\alpha\varphi^{m-k}(x)=x$ for all $x\in\varphi^k(\D_n)$. Let 
$r\geq k$ be so large that $\alpha\varphi^{m-k}(\D_n^k)\subseteq\D_n^r$. Since $\D_n^r=
\D_n^k\varphi^k(\D_n^{r-k})$ and $\alpha\varphi^{m-k}$ acts identically on $\varphi^k(\D_n^{r-k})$, 
it follows that $\alpha\varphi^{m-k}(\D_n^r)\subseteq\D_n^r$. Since the map is injective 
and the space finite dimensional, we have $\alpha\varphi^{m-k}(\D_n^r)=\D_n^r$. As $\alpha
\varphi^{m-k}$ acts identically on $\varphi^r(\D_n)$, it follows that there is a permutation unitary 
$u\in\P_n^r$ such that $\alpha\varphi^{m-k}=\ad(u)|_{\D_n}$. In particular, $\alpha\varphi^{m-k}$ is 
an automorphism of $\D_n$. Thus $m-k=0$ and hence $\alpha=\ad(u)|_{\D_n}$. If $m\leq k$ then 
$\alpha^{-1}\varphi^k=\varphi^m$ and we argue in the same way. This proves one implication. 
The other one is obvious. 
\hfill$\Box$

\begin{lemma}\label{braiding}
If $\alpha\in\aut(\D_n)$ then there exists $\beta\in\aut(\D_n)$ such that $\alpha\varphi=\beta\varphi\alpha$ 
and $\beta(x)=\alpha(x)$ for all $x\in\D_n^1$. 
\end{lemma}
{\em Proof.} Given $x\in\D_n$, there exist unique elements $x_1,\ldots,x_n$ in $\D_n$ such that 
$x=\sum_{j=1}^n P_j\varphi(x_j)$. We define
\begin{equation}\label{beta}
\beta(x) = \alpha\left(\sum_{j=1}^n P_j\varphi(\alpha^{-1}(x_j))\right). 
\end{equation}
It follows that $\beta$ is an automorphism of $\D_n$ satisfying the required braiding property. 
\hfill$\Box$

\medskip
If $\alpha\in\aut(\D_n)$ and $\beta\in\aut(\D_n)$ is defined by formula (\ref{beta}) then we 
call $\beta$ the {\em braiding automorphism} for $\alpha$. 

\medskip
As shown in  \cite[Lemma 5]{S}, if $u\in\P_n$ and $\lambda_u$ 
is an automorphism of $\O_n$ then the restriction of $\lambda_u$ to $\D_n$ belongs 
to $\Gf_n$. This can be further generalized, as follows. 

\begin{theorem}\label{asym}
Let $\alpha \in \aut(\D_n)$. Then the following conditions are equivalent:
\begin{itemize}
\item[(1)] $\alpha^{-1}$ has property (P);
\item[(2)] there exists a permutation $u \in \P_n$ such that $\alpha = \lambda_u|_{\D_n}$.
\end{itemize}
In this case, 
$\alpha$ extends to an endomorphism of $\F_n$ and
both $\alpha$ and $\alpha^{-1}$ are $\tau$-preserving.
% Moreover, $\lambda_u$ is an automorphism of $\O_n$ if and only if $\alpha \in \Gf_n$.
\end{theorem}

{\em Proof.}
(2) $\Rightarrow$ (1): 
% This is basically the same argument of the proof of \cite[Lemma 5]{S}, after noticing that 
% it still works just requiring $\lambda_u|_{\D_n} \in \aut(\D_n)$. 
By assumption, $\alpha$ can be written as 
$\lim_{h\to\infty}{\rm Ad}(u_h)$ (pointwise norm limit)
and therefore $\alpha^{-1} = \lim_{h\to \infty}{\rm Ad}(u_h)^*$. Now,
if $u \in \P_n^r$, one has,
for all $k-1 \geq r$ and all $1 \leq i \leq n$,
$$\begin{aligned}
\alpha^{-1}(\varphi^{k-1}(S_i S_i^*)) 
& =  \lim_{h \to \infty} {\rm Ad} (u_h)^*(\varphi^{k-1}(S_i S_i^*)) \\
& =  \lim_{h \to \infty} \varphi^{h-1}(u^*) \ldots \varphi(u^*)u^* 
   \varphi^{k-1}(S_i S_i^*)u \varphi(u) \ldots \varphi^{h-1}(u) \\
& =  \lim_{h \to \infty} \varphi^{h-1}(u^*) \ldots \varphi^{k-r}(u^*)
   \varphi^{k-r}(\varphi^{r-1}(S_i S_i^*))\varphi^{k-r}(u) \ldots \varphi^{h-1}(u) \\
& =  \varphi^{k-r}\Big(\lim_{h \to \infty}\varphi^{h-k+r-1}(u^*) \ldots u^* 
   \varphi^{r-1}(S_i S_i^*)u \ldots \varphi^{h-k+r-1}(u)\Big) \\
& =  \varphi^{k-r}(\alpha^{-1}(\varphi^{r-1}(S_i S_i^*))) \ .
\end{aligned}$$
\\
(1) $\Rightarrow$ (2): 
At first we observe that there exists a permutation unitary $u\in\P_n$ 
such that the braiding automorphism for $\alpha$ is of the form $\beta=\Ad(u)|_{\D_n}$.
Indeed, let $\alpha^{-1}$ satisfy (P) with $m$. Since $\alpha^{-1}\varphi^{m+1}=
\varphi \alpha^{-1}\varphi^m$, we have $\varphi^{m+1} = \alpha \varphi \alpha^{-1} \varphi^m
= \beta \varphi^{m+1}$. Therefore, there exists a permutation unitary $u\in\P_n$ such that 
$\beta=\Ad(u)|_{\D_n}$, by Lemma \ref{amk}.

Now we verify that the unitary $u$ as above satisfies $\alpha=\lambda_u|_{\D_n}$. Note that, 
by induction, $\alpha \varphi^k = (\beta\varphi)^k \alpha$ for all $k$. Therefore,
$\alpha\varphi^k = (\Ad(u)\varphi)^k\alpha = \Ad(u_k)\varphi^k\alpha$ for all $k$. 
In particular, $\alpha\varphi^k(P_i) = \Ad(u_k)\varphi^k\alpha(P_i) = \Ad(u_{k+1})\varphi^k(P_i)$ 
for all $k$ and $i$.
Now we compute
$$\begin{aligned}
\alpha(P_{i_1}\varphi(P_{i_2}) \cdots \varphi^{k-1}(P_{i_k}))
& = \alpha(P_{i_1}) \alpha\varphi(P_{i_2}) \cdots \alpha \varphi^{k-1}(P_{i_k}) \\
& = \Ad(u)(P_{i_1}) \Ad(u_2)\varphi(P_{i_2}) \cdots \Ad(u_k)\varphi^{k-1}(P_{i_k}) \\
& = u_k P_{i_1} u_k^* u_k \varphi(P_{i_2})u_k^* \cdots u_k\varphi^{k-1}(P_{i_k})u_k^* \\
& = \Ad(u_k) (P_{i_1}\varphi(P_{i_2}) \cdots \varphi^{k-1}(P_{i_k})). 
\end{aligned}$$
That is, $\alpha = \lambda_u|_{\D_n}$, as required.

If the above conditions are satisfied, then in particular $\alpha$ extends to an endomorphism 
$\lambda_u|_{\F_n}$ of $\F_n$, Thus $\tau\lambda_u|_{\F_n}=\tau$ by uniqueness of trace 
on $\F_n$ (or the fact that $\lambda_u|_{\F_n}$ is a point-wise limit of inner automorphisms). 
Hence $\tau \alpha = \tau$ as well. This clearly implies also that $\tau \alpha^{-1}=\tau$.
% \\
% If $\alpha \in \Gf_n$ by the first part of the proof there exists also
% a permutation $v$ such that $\alpha^{-1}=\lambda_v|_{\D_n}$, and thus
% $\lambda_u \lambda_v \in \aut(\O_n)$ by \cite[Proposition 3.2]{Co}.
\hfill $\Box$

\medskip
The preceding theorem provides a dynamical explanation 
of the puzzling phenomenon observed in \cite{CS} through combinatorial arguments,
namely the coexistence of permutative automorphisms and proper endomorphisms of $\O_n$ 
restricting to automorphisms of the diagonal,
in the form of a different dynamics they induce on $X_n$.

\medskip
The following proposition is yet another slight generalization of \cite[Lemma 5]{S}. We omit the proof.

\begin{proposition}\label{Ploc}
Let $\lambda_u$ be a localized automorphism of $\O_n$. 
Then, there exists some nonnegative integer $m$ such that, for all $k \geq m$ and 
all $x\in\F_n^1$, 
\begin{equation}
\lambda_u^{-1} \circ \varphi^k(x) = \varphi^{k-m} \circ \lambda_u^{-1} \circ \varphi^m(x) \ .
\end{equation}
\end{proposition}

\begin{theorem}\label{extendingautos}
If $\alpha\in\Gf_n$ then there exists a permutation unitary $u\in\P_n$ such that 
$\alpha=\lambda_u|_{\D_n}$. If $w$ is any unitary in $\O_n$ such that $\lambda_w|_{\D_n}=\alpha$, 
then $\lambda_w\in\aut(\O_n)$. 
\end{theorem}
{\em Proof.} The first statement follows immediately from Theorem \ref{asym}.
Let $w\in\U(\O_n)$ be such that $\lambda_w|_{\D_n}=\alpha$. Then, by the same argument, 
there exists a permutation unitary $v\in\P_n$ such that $\lambda_v|_{\D_n}=\alpha^{-1}$. Thus  
$\lambda_w\lambda_v$ is an endomorphism of $\O_n$ acting identically on $\D_n$, and 
consequently $\lambda_w\in\aut(\O_n)$ by \cite[Proposition 3.2]{Co}. 
\hfill$\Box$

\begin{corollary}\label{uhf}
Each $\alpha\in\Gf_n$ can be extended to an automorphism of $\F_n$. 
% In particular, $\tau\alpha=\tau$. 
\end{corollary}

\medskip
By Theorem \ref{asym}, the restriction $r:\aut(\O_n,\D_n)\to\aut(\D_n)$ 
yields a group embedding $\lambda(\P_n)^{-1}\hookrightarrow\Gf_n$, see \cite{Cun2,CS}. 
Since the restriction map $r:\aut(\O_n,\D_n)\to\aut(\D_n)$ is injective on $\lambda(\P_n)^{-1}$, 
\cite{Cun2,CS}, Theorem \ref{extendingautos} yields the following. 

\begin{corollary}\label{permautos}
The restriction $r:\lambda(\P_n)^{-1}\to\Gf_n$ is a group isomorphism. 
\end{corollary}

We recall from \cite{Cun2} that $\aut(\O_n,\D_n)$ is the normalizer of $\aut_{\D_n}(\O_n)$ in $\aut(\O_n)$ and
it can be also described as the group $\lambda(\N_{\O_n}(\D_n))^{-1}$ of automorphisms of $\O_n$ induced by elements in the (unitary) normalizer $\N_{\O_n}(\D_n)$.
Furthermore, using \cite{Pow}, 
one can show that $\aut(\O_n,\D_n)$ has the structure of a semidirect product $\aut_{\D_n}(\O_n) \rtimes \lambda(\SS_n)^{-1}$ \cite{CS}.
In particular, the group $\lambda(\P_n)^{-1}$ is isomorphic with the quotient of 
the group $\aut(\O_n,\D_n)\cap\aut(\O_n,\F_n)$ by its normal subgroup $\aut_{\D_n}(\O_n)$. 
We call it the {\em restricted Weyl group} of $\O_n$, cf. \cite{Cun2,CS}. Thus, the preceding 
corollary provides in a sense an answer to the question raised by Cuntz in \cite{Cun2}. 

\begin{corollary}\label{sngn}
Let $u\in\SS_n$ be such that $\lambda_u(\D_n)=\D_n$ and $\lambda_u|_{\D_n}\in\Gf_n$. Then $u\in\P_n$. 
\end{corollary}
{\em Proof.} By Theorem \ref{extendingautos}, there is $w\in\P_n$ such that $\lambda_w(x)=
\lambda_u(x)$ for all $x\in\D_n$. However, the restriction map $r:\lambda(S_n)^{-1}\to\aut(\D_n)$ is 
injective, \cite{Cun2,CS}. Thus $u=w$ belongs to $\P_n$. 
\hfill$\Box$

%%%%%%%%%%%%%%%%%%%%%%%%%%%%%%%%%%%%%%%%%%%%%%
%%%%%%%%%%%%%%%%%%%%%%%%%%%%%%%%%%%%%%%%%%%%%%

\section{The restricted outer Weyl group of the Cuntz algebra and shift endomorphisms}

Denote by $\operatorname{Inn}\lambda(\P_n)^{-1}$ the normal subgroup of $\lambda(\P_n)^{-1}$ 
consisting of all inner permutative automorphisms $\{\Ad(u):u\in\P_n\}$. We call the quotient 
$\lambda(\P_n)^{-1}/\operatorname{Inn}\lambda(\P_n)^{-1}$ the {\em restricted outer Weyl group} 
of $\O_n$. From Theorem \ref{extendingautos}, we get the following. 

\begin{corollary}\label{weyl}
The restricted outer Weyl group of $\O_n$ is naturally isomorphic to the quotient $\Gf_n/\iGf_n$. 
\end{corollary}
In what follows, if $\alpha\in\Gf_n$ then we denote its class in $\Gf_n/\iGf_n$ by $\overline{\alpha}$.
We denote by $\en(\D_n,\varphi)$ the semigroup of unital, injective $*$-homomorphisms from $\D_n$ 
into itself which commute with the shift. We define $E_n$ as the collection of all those 
$\alpha\in\en(\D_n,\varphi)$ for which there exists an $m$ and a $\beta\in\en(\D_n,\varphi)$  such that 
\begin{equation}\label{phiinverse}
\alpha\beta=\varphi^m. 
\end{equation}
In such a case, we have $\alpha\beta\alpha=\varphi^m\alpha=\alpha\varphi^m$, and thus injectivity 
of $\alpha$ implies that $\beta\alpha=\varphi^m$ as well. In particular, $\beta$ itself belongs to $E_n$. 
$E_n$ is a subsemigroup of $\en(\D_n,\varphi)$ containing all powers of $\varphi$, $\Gf^0$, 
as well as all endomorphisms $\alpha\varphi^m$ with $\alpha\in\Gf_n$ and suitably large $m$. 

We note that if $\alpha\in E_n$ then $\alpha_*:X_n\to X_n$ is an open mapping. Indeed, if $U\subseteq 
X_n$ is open then so is $V=\beta_*^{-1}(U)$, and $U=\beta_*(V)$ since $\beta_*$ is surjective. 
Then $\alpha_*(U)=\alpha_*\beta_*(V)=\varphi_*^m(V)$ is open, since $\varphi_*^m$ is an 
open mapping. 

\begin{lemma}\label{shifttrace}
If $\alpha\in E_n$ then $\alpha_*$ acts bijectively on periodic words. Consequently, 
$\tau\alpha(d)=\tau(d)$ for all $d\in\D_n$. 
\end{lemma}
{\em Proof.} Say $x\in X_n$ has period $r$ if $\varphi_*^r(x)=x$. For each $r$, the set $X_n(r)$ of 
all points with period $r$ is finite. Furthermore, $\varphi_*$ restricts to a {\em bijection} on each $X_n(r)$. 
If $\mu\in W_n$ and $r\geq|\mu|$ then there are exactly $n^{r-|\mu|}$ words $x\in X_n(r)$  
such that $P_\mu(x)=1$. 

Now let $\alpha$ be in $E_n$, and let $\beta\in E_n$ be such that (\ref{phiinverse}) holds. 
For each $r$ we have  $\alpha_*(X_n(r))\subseteq X_n(r)$, since $\alpha_*$ commutes with 
$\varphi_*$. If $x,y\in X_n(r)$ and $\alpha_*(x)=\alpha_*(y)$ then also 
$\varphi_*^m(x)=\beta_*\alpha_*(x)=\beta_*\alpha_*(y)=\varphi_*^m(y)$. Thus 
$x=y$, since $\varphi_*$ acts bijectively on $X_n(r)$. Therefore, $\alpha_*$ yields a one-to-one 
mapping from $X_n(r)$ to itself. By finiteness of $X_n(r)$, this map is bijective. 

Now let $\mu\in W_n$ and let $\alpha(P_\mu)=\sum_{j=1}^k P_{\mu_j}$. By subdividing, if necessary, 
we may assume that each $\mu_j$ is of the same length $r$ and that $r\geq|\mu|$. Then the number 
of words $x\in X_n(r)$ such that $P_\mu(x)=1$ is the same as the number of words $y\in X_n(r)$ 
such that $P_{\mu_j}(y)=1$ for some $j$. Thus $k=n^{r-|\mu|}$, and consequently 
$$ \tau\alpha(P_\mu)=\frac{k}{n^r}=\frac{1}{n^{|\mu|}}=\tau(P_\mu). $$
Since $\mu$ was arbitrary, the proposition follows. 
\hfill$\Box$

\medskip
The following lemma is due to Mike Boyle, \cite{Boy}, although the proof given below is our own. 

\begin{lemma}[M. Boyle]\label{boyle}
If $\alpha\in E_n$ then there exists a $k$ such that the mapping $\alpha_*$ is  $k$-to-one. This 
$k$ divides a power of $n$. Thus if $n$ is prime then there exists an $r$ such that 
$\alpha_*$ is $n^r$-to-one. 
\end{lemma}
{\em Proof.} Let $\alpha,\beta\in E_n$ satisfy (\ref{phiinverse}). Since $\beta_*\alpha_*=
\varphi_*^m$ is $n^m$-to-one, each $y\in X_n$ has at most $n^m$ 
inverse images under $\alpha_*$. Let $k$ be the minimal cardinality of $\alpha_*^{-1}(y)$, and 
let $\omega\in X_n$ be such that $\alpha_*^{-1}$ has exactly $k$ elements. Then the set 
$\alpha_*^{-1}(\varphi_*^{-1}(\omega))=\varphi_*^{-1}(\alpha_*^{-1}(\omega)$ has $nk$ 
elements. By the minimality of $k$, this can only happen if each element of $\varphi_*^{-1}(\omega)$ 
has $k$ inverse images under $\alpha_*$. Let $\Omega$ be the smallest subset of $X_n$ 
containing $\omega$ and closed under taking inverse images of $\varphi_*$. It follows from the above 
that for each $y\in\Omega$ the set $\alpha_*^{-1}(y)$ has $k$ elements. Clearly, $\Omega$ is 
dense in $X_n$. Now let $y\in X_n$ be arbitrary, and let $\alpha_*^{-1}(y)$ have $k'$ elements. 
Since $\alpha_*$ is an open mapping, there exists an open subset $V$ of $X_n$ containing $y$ whose 
each element has at least $k'$ inverse images under $\alpha_*$. Since $V\cap\Omega\neq\emptyset$,  
we have $k'=k$. 

Applying the same reasoning to $\beta$ instead of $\alpha$ we get an $l$ such that $\beta_*$ is 
$l$-to-one. Thus $kl=n^m$. 
\hfill$\Box$

\medskip
We define an equivalence relation $\sim$ in $E_n$ as follows: $\alpha\sim\beta$ if there exists a $k$ 
such that either $\alpha=\beta\varphi^k$ or $\alpha\varphi^k=\beta$. Then we set 
$$ \Ef_n = E_n/\sim. $$
By construction, $\Ef_n$ is a group. If $\alpha\in E_n$ then we denote its class in $\Ef_n$ by 
$\overline{\alpha}$. Let $\aut(\Sigma_n)$ denote the group of automorphisms of the full two-sided 
$n$-shift, and let $\la \sigma \ra$ be its subgroup generated by the two-sided shift $\sigma$. 
The following proposition is well-known, but for completeness we include a proof. 

\begin{proposition}\label{quotientiso}
The groups $\Ef_n$ and $\aut(\Sigma_n)/\la \sigma \ra$ are isomorphic. 
\end{proposition}
{\em Proof.}
We can realize $\D_n$ as $\otimes_{i \in {\mathbb N}} \D_n^1$.
Likewise, we consider $\widetilde\D_n := \otimes_{i \in {\mathbb Z}} \D_n^1$,
equipped with the two-sided shift automorphism $\widetilde\varphi$.
The Gelfand spectrum of $\widetilde{\D}_n$ can be identified with $\Sigma_n$ and
the canonical embedding $\D_n \to \widetilde{\D}_n$ corresponds to the canonical surjection 
$\Sigma_n \to X_n$.

Now, every shift invariant endomorphism $\rho$ of $\D_n$ canonically extends 
to an endomorphism $\tilde\rho$ of $\widetilde\D_n$, uniquely determined
by the properties of being $\widetilde\varphi$-invariant
and restricting to $\rho$ on $\D_n \subset \widetilde\D_n$.
Of course, the shift endomorphism $\varphi$ extends to $\widetilde\varphi$.
If $\rho$ is injective then $\tilde\rho$ is injective too. Moreover, the map
$\rho \mapsto \tilde\rho$ gives a semigroup homomorphism from 
${\rm End}(\D_n,\varphi)$ to ${\rm End}(\widetilde\D_n,\widetilde\varphi)$, which is clearly injective.

If $\rho \in E_n$ then it is easy to see that $\tilde\rho$ is surjective and thus it is an
automorphism of $\widetilde\D_n$. Therefore, we get an injective semigroup homomorphism
$E_n \to {\rm Aut}(\widetilde\D_n,\widetilde\varphi)$. Passing to quotients, the previous 
map provides a well-defined and injective group homomorphism from $\Ef_n$ into 
${\rm Aut}(\widetilde\D_n,\widetilde\varphi)/\langle \widetilde\varphi \rangle
\cong \aut(\Sigma_n)/\la \sigma \ra$.

In order to show surjectivity of this map, we observe that
given any $\widetilde\varphi$-commuting automorphism $\eta$ of $\widetilde\D_n$
there exists a nonnegative integer $k$ such that $\eta':=\widetilde\varphi^k \eta |_{\D_n}$ 
belongs to $E_n$ (use the fact that $\Sigma_n$ is the Cantor set). Furthermore, thanks to 
uniqueness of the extension, one has $\widetilde\varphi^{-k}\tilde{\eta'}=\eta$ and the 
proof is complete.
\hfill$\Box$

\begin{theorem}\label{endogroup}
There exists an embedding of $\Gf_n/\iGf_n$ into $\aut(\Sigma_n)/\la \sigma \ra$. If $n$ is 
prime then these two groups are isomorphic. 
\end{theorem}
{\em Proof.} By virtue of Proposition \ref{quotientiso}, we may replace $\aut(\Sigma_n)/\la \sigma \ra$ 
with $\Ef_n$. 

 Let $\alpha\in\Gf_n$ satisfy (P) with $m$. Then $\alpha\varphi^m\in E_n$, and we 
map $\alpha$ to $\overline{\alpha\varphi^m}$ in $\Ef_n$. Clearly, this definition does not depend on 
the choice of $m$ and thus the map is well defined. One easily checks that this map is a group homomorphism. 
If $\alpha\in\iGf_n$ then $\alpha\varphi^k=\varphi^k$ for sufficiently large $k$, and thus the image 
of such $\alpha$ in $\Ef_n$ is the trivial element. Thus, the homomorphism $\Gf_n\to\Ef_n$ factors 
though $\Gf_n/\iGf_n$, and we get a homomorphism $\Gf_n/\iGf_n\to\Ef_n$. The latter map is 
injective. Indeed, let $\alpha\in\Gf_n$ satisfy (P) with $m$ and $\overline{\alpha\varphi^m}=
\overline{\id}$. Then there is a $k$ such that either $\alpha\varphi^{m+k}=\id$ or 
$\alpha\varphi^m=\varphi^k$. In either case, $\alpha\in\iGf_n$ by Lemma \ref{amk}. 

Now assuming $n$ prime we show that the map $\Gf_n/\iGf_n\to\Ef_n$ is surjective. Let $\alpha\in E_n$ 
and let $\beta\in E_n$ be such that identity (\ref{phiinverse}) holds. There is an $r$ such 
that $\alpha_*$ is $n^r$-to-one, by Lemma \ref{boyle}. Since $\alpha_*$ is an open map 
onto the Cantor set it admits a continuous section $f_1$, 
\cite[Corollary 1.4]{M}. Space $X_n\setminus f_1(X_n)$ being open is completely metrizable and thus 
the restriction of $\alpha_*$ to $X_n\setminus f_1(X_n)$ admits a section $f_2$. Continuing in this manner, 
we arrive at a maximal set of independent, continuous sections. 
For convinience, we label these sections with words from $W_n^r$.  That is, there exist continuous 
functions $f_\mu:X_n\to X_n$, $\mu\in W_n^r$, such that $\alpha_* f_\mu=\id$ and 
$\alpha_*^{-1}(x)=\{f_\mu(x) : \mu\in W_n^r \}$ for each $x\in X_n$. Define a mapping 
$g:X_n \to X_n$ by $g(\mu x)=f_\mu(x)$. Then $g$ is a homeomorphism of $X_n$. Furthermore,
 we have $\alpha_* g(\mu x) = \alpha_* f_\mu(x) = x = \varphi_*^r(\mu x)$. Thus 
$\alpha_* g = \varphi_*^r$ and hence $\alpha_* = \varphi_*^r g^{-1}$. Let $\psi\in\aut(\D_n)$ 
be such that $\psi_*=g^{-1}$. Then $\alpha = \psi \varphi^r$ and whence $\psi$ satisfies 
condition (P).  Arguing in the same way, we obtain an $\eta\in\aut(\D_n)$ such that $\beta=
\eta\varphi^t$ for some $t$. Then $\varphi^m=\alpha\beta=\psi\eta\varphi^{r+t}$. Thus 
$\psi\eta\in\iGf_n$ by Lemma \ref{amk}. Now $\psi^{-1}$ being a composition of $\eta$ and 
an element from $\iGf_n$, itself satisfies condition (P). Consequently, $\psi$ belongs to $\Gf_n$. 
Clearly, the homomorphism $\Gf_n/\iGf_n\to\Ef_n$ maps $\overline{\psi}$ to $\overline{\alpha}$, 
and the proof is complete. 
\hfill$\Box$

\begin{remarks}
\rm In view of the preceding theorem, the restricted outer Weyl group of $\O_n$ 
has a number of striking properties known to hold for $\aut(\Sigma_n)/\la \sigma \ra$, see \cite{K,LM} 
and the references therein. This is immediate at least for $n$ prime.  For example, we now know 
that in the case of $n=2$ the group $\Gf_n/\iGf_n$ is non-amenable (for $n\geq3$ this has been 
already observed in \cite{S,CS}).  We wonder if the extensive theory of shift automorphisms could not 
bring new insight into the structural properties of not just Cuntz algebras but also graph algebras and 
possibly even a larger class of $C^*$-algebras. In particular, it appears to be an intriguing possibility of 
translating some features of the beautiful analysis of the action of automorphisms on periodic points, 
\cite{BK}, into a more algebraic setting.  

In general, when $n$ is not prime, the embedding from Theorem \ref{endogroup} 
is not surjective. This is due to existence of factorizations $\varphi=\alpha\beta$ with $\alpha,\beta\in 
E_n$ and neither $\alpha$ nor $\beta$ being an automorphism, \cite{BM}. Then it is easy 
to verify that $\overline{\alpha}$, $\overline{\beta}$ do not belong to the range of the embedding. 
Nevertheless, even for $n$ not a prime, Theorem \ref{endogroup} sheds a lot of light on the structure of the 
restricted outer Weyl group of $\O_n$.  In particular, it implies that $\Gf_n/\iGf_n$ is residually finite. 
As an example, we give an elementary, self-contained proof of the fact that $\Ef_n$ is residually finite in 
Proposition \ref{rf} below.
\end{remarks}

\begin{proposition}\label{rf}
The group $\Ef_n$ is residually finite. 
\end{proposition}
{\em Proof.} Let $\per$ denote the set of periodic points in $X_n$. 
If $\alpha\in E_n$ then $\alpha_*$ acts bijectively on $\per$, 
Lemma \ref{shifttrace}. Thus the group $\Ef_n$ acts on the orbits of $\per$ 
under the action of $\varphi_*$. The restriction of this action to the orbits contained 
in $X_n(r)$ yields a homomorphism from $\Ef_n$ into a finite permutation group. 
Thus it suffices to show that if $\alpha\in E_n$ is not a power of the shift then  
$\alpha_*$ moves at least one orbit. 

Let $\alpha\in E_n$. Suppose that $\alpha_*$ fixes every orbit of periodic 
points under the action of the shift. This means that for each periodic point $x$  
there exists $k$ such that $\alpha_*(x)=\varphi_*^k(x)$. We claim that 
$\alpha$ is a power of the shift. The proof involves the following five steps. 

1. There is a map $a:W_n^r \to W_n^1$ such that if $y=\alpha_*(x)$ then 
$y_1=a(x_1,\ldots,x_r)$, \cite{H}. 

2. If $x_1=\ldots = x_r = j$ then $a(j,\ldots,j)=j$, for otherwise 
$(j,j,\ldots)$ would be moved by $\alpha_*$, being a fixed point for 
$\varphi_*$. 

3. Let $Z$ be the set of all periodic points of the following form:
$$ (x_1,\ldots,x_t,w,\ldots,w)(x_1,\ldots,x_t,w,\ldots,w)\ldots, $$ 
where $t>2r$, the number of 
$w$'s (in one block) is $d>t$, and $w\neq a(x_1,\ldots,x_r)$. Note that $Z$ 
is dense in $X_n$. By hypothesis on $\alpha_*$, the image under $\alpha_*$ of 
such a periodic point from $Z$ is of the form: 
$$ \begin{aligned} 
\text{either}\;\;\; & (w,\ldots,w,x_1,\ldots,x_t,w,\ldots,w) (w,\ldots,w,x_1,\ldots,x_t,w,\ldots,w)\ldots \\  
\text{or}\;\;\; & (x_p,\ldots,x_t,w,\ldots,x_1,\ldots,x_{p-1}) (x_p,\ldots,x_t,w,\ldots,x_1,\ldots,x_{p-1})\ldots
\end{aligned} $$
However, the former is impossible since $w\neq a(x_1,\ldots,x_r)$ and $t>r$. 
Furthermore, there is a $k$ in $\{1,\ldots,r\}$ such that $x_k\neq w$, for 
otherwise $a(x_1,\ldots,x_r)=w$, contrary to the assumption. Consequently, 
$1 \leq p \leq 2r-1$. 

4. By 3. above, for each $x\in Z$ there is a $0 \leq p \leq 2r-2$ such that 
$\alpha_*(x)
= \varphi_*^p(x)$. Now let $x$ be an arbitrary point in $X_n$. Take a sequence 
$x_k\in Z$ converging to $x$. Passing to a subsequence we may assume that 
there is a fixed $p \leq 2r-2$ such that for each $k$ we have $\alpha_*(x_k)=
\varphi_*^p(x_k)$. By continuity, $\alpha_*(x)=\varphi_*^p(x)$. 

5. For $j=0,1,\ldots,2r-2$ let $Y_j=\{ x\in X_n : \alpha_*(x)=\varphi_*^j(x) \}$. 
By 4. above, $X_n = Y_0 \cup \ldots \cup Y_{2r-2}$, and each $Y_j$ is a closed set. 
Thus there is $j$ such that $Y_j$ contains a non-empty open subset.  
Then there exists a finite word $\mu$ such that 
for all infinite words $x$ we have $\mu x \in Y_j$ and  $\alpha_*(\mu x)=
\varphi_*^j(\mu x)$. Since $\alpha_*$ commutes with $\varphi_*$, we have 
$$ \alpha_*(x) = \alpha_*\varphi_*^{|\mu|}(\mu x) = \varphi_*^{|\mu|}\alpha_*(\mu x) = 
\varphi_*^{j+|\mu|}(\mu x) =\varphi_*^j(x). $$ 
Since $x$ was arbitrary, $\alpha_* = \varphi_*^j$. 
\hfill$\Box$

%%%%%%%%%%%%%%%%%%%%%%%%%%%%%%%%%%%%%%%%%%%%%
%%%%%%%%%%%%%%%%%%%%%%%%%%%%%%%%%%%%%%%%%%%%%

\section{Shift automorphisms}

%%%%%%%%%%%%%%%%%%%%%%%%%%%%%%%%%%%%%%%%%%%%%
%%%%%%%%%%%%%%%%%%%%%%%%%%%%%%%%%%%%%%%%%%%%%

Note that if $\alpha\in\aut(\D_n)$ satisfies (P) with $m=0$ then $\alpha_*$ is just an automorphism 
of the full one-sided $n$-shift (see \cite{K}). The collection of all such automorphisms constitutes 
a subgroup of $\Gf_n$, which we denote $\Gf_n^0$. In the case of $n=2$ we have $\Gf_2^0\cong\Zb_2$ 
(generated by the restriction of Archbold's flip-flop, \cite{Ar}) by \cite[Theorem 3.1.1]{K}. But 
for $n\geq3$ the group $\Gf_n^0$ is infinite (see \cite[Chapter 3]{K} and \cite{BFK,H}). 

\begin{example}\label{Kitchensexample}
\rm Consider an order two automorphism $\alpha$ of $\D_3$ such that $\alpha_*$ changes 
subwords $13$ and $23$ (of any one-sided infinite word) into $23$ and $13$, respectively, as 
in \cite[Example 3.3.10]{K}. Then $\alpha$ belongs to $\Gf_3^0$. Define 
$$ u=P_{11}+P_{12}+P_{21}+P_{22}+P_3+S_{23}S_{13}^*+S_{13}S_{23}^*, $$ 
an order two permutation unitary in $\P_3^2$. We have 
\begin{eqnarray*}
uP_1u^* & = & P_{11}+P_{12}+P_{23}, \\
uP_2u^* & = & P_{21}+P_{22}+P_{13}, \\
uP_3u^* & = & P_3. 
\end{eqnarray*}
One checks that the unitary $u$ commutes with the three minimal projections in $\varphi(u\D_3^1u^*)$. 
This implies (via an easy inductive argument) that the restriction of $\lambda_u$ to $\D_3$ 
commutes with the shift $\varphi$. It follows that $\alpha=\lambda_u|_{\D_3}$. 
\hfill$\Box$\end{example}

We note, in passing, that among the three rooted trees associated with automorphism $\lambda_u$ 
from Example \ref{Kitchensexample} as in \cite[Section 4.1]{CS}, two are of height 2 and one is 
of height 1 (cf. \cite[Section 2.1]{CKS}). 

\medskip
By Theorem \ref{extendingautos}, we already know that each $\alpha\in\Gf_n^0$ may be extended 
to an automorphism of $\O_n$. Below, we provide an alternative proof of this fact, 
involving a thourough description of the underlying structure in this specific case and thus leading to
a more explicit construction of the required permutation unitary. We will need the following two lemmas. 
A straightforward proof of the former is omitted. 

\begin{lemma}\label{shiftext}
Let $\alpha \in \Gf_n^0$.
If $u$ is a unitary in $\O_n$ such that $u x u^* = \alpha(x)$ for all $x \in \D_n^1$
and $u \varphi^j(\alpha(x)) u^* = \varphi^j(\alpha(x))$ for all $x \in \D_n^1$ and all $j=1,\ldots,k-1$
then $u_k x u_k^* = \alpha(x)$ for all $x \in \D_n^k$.
\end{lemma}

\begin{lemma}\label{phicommute}
Each $\alpha\in\Gf_n^0$ commutes with the left-inverse $\phi$ of $\varphi$.
\end{lemma}
{\em Proof.} Choose $r$ so large that $\alpha(\D_n^1) \subseteq \D_n^r$. 
We write $\alpha(P_j) = \sum_{\mu \in W(j)} P_\mu$, where $W(j) \subseteq W_n^r$ is 
a subset of cardinality $n^{r-1}$, since $\tau\alpha=\tau$ by Lemma \ref{shifttrace}.
Clearly, the sets $\{W(j):j=1,\ldots,n\}$ form a partition of $W_n^r$.

Since $\alpha$ commutes with the shift $\varphi$, in order to show that $\alpha$ commutes with 
its left inverse $\phi$ it suffices to prove that $1=\sum_{i=1}^n S_i^* \alpha(P_j) S_i$ for 
all $j=1,\ldots,n$. This is equivalent to $\{(\mu_2,\ldots,\mu_r) :  \mu \in W(j)\} = W_n^{r-1}$. 
Now the last claim will follow if we can show that
$\alpha(P_j) \varphi(P_{i_1} \cdots \varphi^{r-2}(P_{i_{r-1}})) \neq 0$
for all $i_1,\ldots,i_{r-1}$ and all $j$. But the last expression is different from $0$ 
if and only if
$P_j \varphi(\alpha^{-1}(P_{i_1} \cdots \varphi^{r-2}(P_{i_{r-1}})))$ is, and 
this is clearly the case. 
\hfill$\Box$

\begin{theorem}\label{shiftautos}
Let $\alpha$ be an automorphism of the full one-sided $n$-shift, i.e. $\alpha\in\Gf_n^0$. 
Then there exists a permutation unitary $u\in\P_n$ such that $\alpha=\lambda_u|_{\D_n}$.
\end{theorem}
{\em Proof.} Let $\alpha(P_j) = \sum_{\mu \in W(j)} P_\mu$, as in Lemma \ref{phicommute}. 
Then $\{(\mu_2,\ldots,\mu_r) \ | \ \mu \in W(j)\}=W_n^{r-1}$, for every $j=1,\ldots,n$.
Therefore, there exists a unique $u \in \P_n^r$ such that 
$u P_{\mu} u^* = P_{\mu(j)}$ for all $\mu \in W_n^r$,
where $\mu_1 = j$ and
$\mu(j)$ is the unique multi-index in $W(j)$ such that 
$$(\mu_2,\ldots,\mu_r)=(\mu(j)_2,\ldots,\mu(j)_r) \ . $$
It is straightforward to check that the conditions of Lemma \ref{shiftext} are satisfied 
(for arbitrarily large $k$), and hence $\lambda_u|_{\D_n}=\alpha$. 
\hfill$\Box$

\begin{remark}
\rm All the permutation unitaries $u \in \P_n^r$ as in the preceeding theorem have the 
following general structure: there are $n$ functions $f_1, \ldots, f_n$ from 
$W_n^{r-1}$ into $W_n^1$ such that for each $w \in W_n^{r-1}$ 
the map $j \mapsto f_j(w)$ is a permutation of $W_n^1$, and moreover 
for all $\mu \in W_n^r$, $u P_\mu u^*=P_{\sigma_f(\mu)}$, where
$\sigma_f((\mu_1,\ldots,\mu_r))=(f_{\mu_1}(\mu_2,\ldots,\mu_r),\mu_2,\ldots,\mu_r)$.
Therefore,
$$u = \sum_{\mu\in W_n^r} S_{\sigma_f(\mu)}S_\mu^* \ . $$
\end{remark}

\begin{proposition}\label{outershiftauto}
If $\alpha\in\aut(X_n)$ is a non-trivial automorphism of the full one-sided $n$-shift and $\tilde{\alpha}$ is an 
automorphism of $\O_n$ extending $\alpha$, then $\tilde{\alpha}$ is outer. 
\end{proposition}
{\em Proof.} Suppose $\tilde{\alpha}=\Ad(u)$ for some $u\in\U(\O_n)$. Then $u\in\N_{\O_n}(\D_n)$ and 
thus $u=vw$, where $v\in\SS_n$ and $w\in\U(\D_n)$, \cite{Pow}. But $\Ad(w)$ acts identically 
on $\D_n$. Thus $\tilde{\alpha}=\Ad(v)$. Since $\tau\alpha=\tau$, we must have $v\in\P_n^k$ for 
some $k$. But then $\varphi^k(x)=\Ad(v)\varphi^k(x)=\varphi^k\Ad(v)(x)$ for all $x\in\D_n$. Thus 
$v=1$ and consequently $\tilde{\alpha}=\id$. 
\hfill$\Box$

\medskip
We have shown, above, that every automorphism of $\D_n$ commuting with the shift extends to an 
automorphism of $\O_n$. In this subsection we observe that in most cases such an extension cannot 
commute with the shift $\varphi$. 

\begin{proposition}
Let $\lambda_v$ be an endomorphism of $\O_n$. Then $\lambda_v$ commutes with $\varphi$ 
on $\O_n$ if and only if
\begin{equation}
v \varphi(v) \theta \varphi(v^*) = \varphi(v) \theta \ , 
\end{equation}
where $\theta$ is the flip unitary. If moreover $\lambda_v \in \aut(\O_n)$ then
$\lambda_v \varphi = \varphi \lambda_v$ if and only if $v\in\U(\F_n^1)$, i.e. $\lambda_v$ is 
a Bogolubov automorphism of $\O_n$. 
\end{proposition}
{\em Proof.}
The first statement follows easily from the composition rule of endomorphisms and the fact that
$\varphi = \lambda_\theta$, where $\theta \in \F_n^2$.
% For u=\theta of course this must be true and reduces to the fact that \theta
% satisfies the usual YBE.

Since $\lambda_v \varphi = {\rm Ad}(v) \varphi \lambda_v$ for all unitaries $v \in \U(\O_n)$,
the assumption implies that ${\rm Ad}(v) \varphi \lambda_v(x) = \varphi \lambda_v(x)$ 
for all $x \in \O_n$ and thus, $\lambda_v$ being an automorphism, 
$v \in \varphi(\O_n)' \cap \O_n = \F_n^1$.
\hfill $\Box$

\begin{proposition}
Let $u \in \U(\O_n)$ and suppose that $\lambda_u(\F_n)=\F_n$. 
If $\lambda_u$ and $\varphi$ commute on $\F_n$
then $\lambda_u$ is a Bogolubov automorphism.
\end{proposition}
{\em Proof.}
Observe that indeed $u \in \F_n$ \cite{CRS}.
By a similar argument as in the previous proposition, 
we get ${\rm Ad}(u) \varphi \lambda_u = \varphi \lambda_u$ on $\F_n$
and therefore $u \in \varphi(\F_n)' \cap \F_n = \F_n^1$.
\hfill $\Box$

\medskip
In the case $n=2$, if $u$ is a permutation for which
% $\lambda_u(\D_2)=\D_2$ and 
$\lambda_u \in \aut(\O_2)$ commutes with $\varphi$ on $\D_2$ 
then $u$ must be the flip automorphism, i.e. Bogolubov permutation.
(And then, in turn, $\lambda_u$ commutes with $\varphi$ on the whole of $\O_2$.)
% can try to formulate it without reference to permutations, however notice
% that the extension from D_n to O_n is not unique 
% -- two such extension differing by an element in the ``maximal torus''
However Example \ref{Kitchensexample} illustrates that this is not true anymore for $\D_n$ with $n>2$, 
i.e. there are permutation automorphisms $\lambda_u$
(which therefore satisfy automatically $\lambda_u(\D_n)=\D_n$) 
commuting with $\varphi$ on $\D_n$ which are not Bogolubov
(and thus, by the above proposition, they do not commute with $\varphi$ on $\F_n$,
let alone on $\O_n$).

%%%%%%%%%%%%%%%%%%%%%%%%%%%%%%%%%%%%%%%%%%%%%%%%
%%%%%%%%%%%%%%%%%%%%%%%%%%%%%%%%%%%%%%%%%%%%%%%%

\smallskip\noindent
Roberto Conti \\
Dipartimenti di Scienze \\
Universit{\`a} di Chieti-Pescara `G. D'Annunzio' \\
Viale Pindaro 42, I--65127 Pescara, Italy \\
E-mail: conti@sci.unich.it \\

\smallskip\noindent
Jeong Hee Hong \\
Department of Applied Mathematics \\
Korea Maritime University \\
Busan 606--791, South Korea \\
E-mail: hongjh@hhu.ac.kr \\

\smallskip \noindent
Wojciech Szyma{\'n}ski\\
Department of Mathematics and Computer Science \\
The University of Southern Denmark \\
Campusvej 55, DK-5230 Odense M, Denmark \\
E-mail: szymanski@imada.sdu.dk

\end{document}